\documentclass[a4paper]{amsart}
\usepackage{epsfig,a4wide}

\newtheorem{theorem}{Theorem}[section]
\newtheorem{proposition}[theorem]{Proposition}
\newtheorem{lemma}[theorem]{Lemma}

\newtheorem{claim}[theorem]{Claim}

\theoremstyle{remark}

\begin{document}

\title{Reducing Dehn fillings and small surfaces}

\author[S. Lee]{Sangyop Lee}
\address{School of Mathematics, Korea Institute for Advanced Study,
         207-43 Cheongryangri-dong, Dongdaemun-gu, Seoul 130-012, Korea}
\email{slee@kias.re.kr}
\author[S. Oh]{Seungsang Oh$^1$}
\address{Department of Mathematics, Korea University,
         1, Anam-dong, Sungbuk-ku, Seoul 136-701, Korea}
\email{soh@math.korea.ac.kr}
\author[M. Teragaito]{Masakazu Teragaito$^2$}
\address{Department of Mathematics and Mathematics Education,
         Hiroshima University, Kagamiyama 1-1-1,
         Higashi-Hiroshima 739-8524, Japan}
\email{teragai@hiroshima-u.ac.jp}
\thanks{$^1$Supported by KOSEF research project No. R05-2001-000-00015-0}
\thanks{$^2$Partially supported by Japan Society for the Promotion of Science,
        Grant-in-Aid for Scientific Research (C) 14540082, 2002.}

\begin{abstract}
In this paper we investigate the distances between Dehn fillings
on a hyperbolic 3-manifold that yield $3$-manifolds containing
essential small surfaces including non-orientable surfaces.
Especially we study the situations
where one filling creates an essential sphere or projective plane,
and the other creates an essential sphere, projective plane,
annulus, M\"{o}bius band, torus or Klein bottle,
all 11 pairs of such non-hyperbolic manifolds.
\end{abstract}

\maketitle

\section{Introduction} \label{sec:intro}

Let $M$ be a compact, connected, orientable $3$-manifold
with a torus boundary component $\partial_0 M$.
Let $\gamma$ be a {\it slope\/} on $\partial_0 M$, that is,
the isotopy class of an essential simple closed curve on $\partial_0 M$.
The $3$-manifold obtained from $M$ by {\em $\gamma$-Dehn filling\/}
is defined to be $M(\gamma)=M\cup V_{\gamma}$,
where $V_{\gamma}$ is a solid torus glued to $M$ along $\partial_0M$
in such a way that $\gamma$ bounds a meridian disk in $V_{\gamma}$.

By a {\it small surface\/} we mean one with non-negative Euler characteristic
including non-orientable surfaces.
Such surfaces play a special role in the theory of $3$-dimensional manifolds.
Thurston's geometrization theorem for Haken manifolds \cite{Th} asserts that
a hyperbolic $3$-manifold $M$ with non-empty boundary
contains no essential small surfaces.
Furthermore, if $M$ is hyperbolic, then the Dehn filling $M(\gamma)$ is also
hyperbolic for all but finitely many slopes \cite{Th},
and a good deal of attention has been directed towards obtaining
a more precise quantification of this statement.

Let us say that a $3$-manifold is of {\em type\/} $\mathcal{S}$,
$\mathcal{D}$, $\mathcal{A}$, or $\mathcal{T}$,
if it contains an essential orientable small surface which is
an essential sphere, disk, annulus or torus,
and of {\em type\/} $\mathcal{P}$, $\mathcal{B}$, or $\mathcal{K}$,
if it contains a non-orientable small surface which is
a projective plane, M\"{o}bius band or Klein bottle, respectively.
Especially $3$-manifolds of type $\mathcal{S}$, $\mathcal{D}$,
$\mathcal{A}$ and $\mathcal{T}$ are called reducible,
$\partial$-reducible, annular and toroidal, respectively.
The distance $\Delta(\gamma_1,\gamma_2)$ between two slopes
on a torus is their minimal geometric intersection number.
The bound $\Delta(X_1,X_2)$ is the least nonnegative number $m$
such that if $M$ is a hyperbolic manifold which admits two
Dehn fillings $M(\gamma_1)$, $M(\gamma_2)$ of type $X_1$, $X_2$,
respectively, then $\Delta(\gamma_1,\gamma_2)\leq m$.
Surveys of the known bounds of various choices $(X_1,X_2)$
and the maximal values realized by known
examples are given in \cite{G2, W3, EW}.

In this paper we consider $M(\gamma_i)$, $i=1,2$ of only the six types
$\mathcal{S}$, $\mathcal{P}$, $\mathcal{A}$, $\mathcal{B}$,
$\mathcal{T}$, or $\mathcal{K}$.
Suppose that $M(\gamma_i)$ contains such a small surface $\widehat{F}_i$.
Then we may assume that $\widehat{F}_i$ meets the attached solid torus
$V_{\gamma_i}$ in a finite collection of meridian disks,
and is chosen so that the number of disks $n_i$ is minimal
among all such surfaces in $M(\gamma_i)$.
The main results of this paper are the followings.

\begin{theorem} \label{thm:SS}
Suppose that $M$ is hyperbolic.
If $M(\gamma_1)$ and $M(\gamma_2)$ are of type $\mathcal{S}$ or $\mathcal{P}$,
then $\Delta(\gamma_1,\gamma_2) \leq 1$.
\end{theorem}

The original proof of the case $(\mathcal{S},\mathcal{S})$ is very complicated \cite{GL1}.
Our proof is remarkably short, although it is based on the analysis of intersections
of two surfaces as well as \cite{GL1}.
Recently, Hoffman and Matignon \cite{HoM} also gives such a short proof in the almost same line, but
ours is still simpler than it.

If a $3$-manifold contains a projective plane, then it is either a reducible manifold or
the $3$-dimensional projective space.
Hence type $\mathcal{P}$ breaks down into type $\mathcal{S}$ and CYC, which means
the class of manifolds with finite cyclic fundamental groups.
References are: \cite{BZ2} for $\Delta(\mathcal{S},\mbox{CYC})=1$; \cite{CGLS} for $\Delta(\mbox{CYC},\mbox{CYC})=1$.
See also \cite{M,Te} for $\Delta(\mathcal{P},\mathcal{P})=1$.

\begin{theorem} \label{thm:SA}
Suppose that $M$ is hyperbolic.
If $M(\gamma_1)$ is of type $\mathcal{S}$ or $\mathcal{P}$,
and $M(\gamma_2)$ is of type $\mathcal{A}$ or $\mathcal{B}$,
then either $\Delta(\gamma_1,\gamma_2) \leq 1$,
or $\Delta(\gamma_1,\gamma_2) = 2$
with $n_2 = 2$ when $M(\gamma_2)$ is of type $\mathcal{A}$
(or $n_2 = 1$ when of type $\mathcal{B}$).
\end{theorem}

For the case $(\mathcal{S},\mathcal{A})$, Wu \cite{W3} showed that $\Delta(\mathcal{S},\mathcal{A})=2$
by using the sutured manifold theory,
and he asked whether $n_2=2$ when $M(\gamma_1)$ is of type $\mathcal{S}$ and
$M(\gamma_2)$ is of type $\mathcal{A}$  with $\Delta(\gamma_1,\gamma_2)=2$ \cite[Question 5.8]{W3}.
Our Theorem \ref{thm:SA} gives the affirmative answer of this question.
Note that type $\mathcal{B}$ breaks down into types $\mathcal{S}$, $\mathcal{A}$ and $\mathcal{D}$.

\begin{theorem} \label{thm:ST}
Suppose that $M$ is hyperbolic.
If $M(\gamma_1)$ is of type $\mathcal{S}$ or $\mathcal{P}$,
and $M(\gamma_2)$ is of type $\mathcal{T}$ or $\mathcal{K}$,
then either $\Delta(\gamma_1,\gamma_2) \leq 2$,
or $\Delta(\gamma_1,\gamma_2) = 3$
with $n_2 = 2$ when $M(\gamma_2)$ is of type $\mathcal{T}$
(or $n_2 = 1$ when of type $\mathcal{K}$).
\end{theorem}

For the case $(\mathcal{S},\mathcal{T})$, Oh \cite{O} and Wu \cite{W2} independently showed that
$\Delta(\mathcal{S},\mathcal{T})=3$.
(See also \cite{LOT} for its short proof.)
Hence Theorem \ref{thm:ST} gives an improvement of this result.
In \cite{JLOT}, we have showed the conclusions for two cases $(\mathcal{P},\mathcal{T})$
and $(\mathcal{P},\mathcal{K})$.

If a $3$-manifold is of type $\mathcal{K}$,
then it is of type $\mathcal{S}$, $\mathcal{T}$ or a Seifert fibered manifold with finite fundamental group (a prism manifold).
Indeed, non-orientable cases are necessary to prove orientable cases in our arguments.
Such phenomenon is also observed in \cite{GL2} and \cite{GL3}.
%%%%%%%%%%%%%%%%%%%%%%%

We should emphasize that all 11 cases can be treated in a unified argument in this paper.
(Also, we have a plan to continue the study for the remaining pairs of $\mathcal{S,P,D,A,B,T,K}$.)
Our main tool in this paper is a two-cornered cycle, which was introduced by Hoffman \cite{H}.
Hoffman showed that the disk bounded by a great $x$-cycle contains a pair of
specific two-cornered cycles, called a seemly pair,
and it can be used to find a new essential sphere and projective plane
meeting the attached solid torus in a fewer times than an original surface, leading to a contradiction.
We define a slight generalization of a great $x$-cycle, called an $x$-face, and
show that it contains a seemly pair by a simpler argument than that in \cite{H} and that
such a pair is useful for the types $\mathcal{A}$ and $\mathcal{B}$ as well.

%%%%%%%%%%%%%
By virtue of Theorem \ref{thm:SS} and the fact $\Delta(\mathcal{S},\mathcal{D})=0$ \cite{S}
(and hence $\Delta(\mathcal{P},\mathcal{D})=0$),
we can put the following assumption throughout the paper to simplify the arguments greatly:

\medskip
\noindent\textbf{Assumption.} If $M(\gamma_i)$ is of type $\mathcal{A},\mathcal{B},\mathcal{T}$ or
$\mathcal{K}$, then we assume that $M(\gamma_i)$ is neither of type $\mathcal{S}$ nor of type $\mathcal{D}$.
\medskip

%%%%%%%%%%%%%%%%%%%%%%%%%%%%%%%%%%%%%%%%%%%%%%%%%%%%%%%%%%%%%%%%%%%%%%%%%%%%%%%%%%%%%%%%%%%%%%%%%%%%%%%
\section{Graphs of surface intersections} \label{sec:graph}

Hereafter $M$ is a hyperbolic $3$-manifold
with a torus boundary component $\partial_0 M$.
An orientable surface properly embedded in $M$ is called \textit{essential\/} if it is
either (i) incompressible, not boundary parallel and not a sphere, or (ii) a sphere
which does not bound a $3$-ball.
Note that any essential small surface is also boundary incompressible by our assumption in Section \ref{sec:intro}.
In this section we describe how (essential) small surfaces
$\widehat{F}_1$ and $\widehat{F}_2$,
in $M(\gamma_1)$ and $M(\gamma_2)$ respectively,
give rise to labelled intersection graphs $G_i \subset \widehat{F}_i$
for $i=1,2$ in a general context.

As in Section \ref{sec:intro},
let $\widehat{F}_i$ be a small surface in $M(\gamma_i)$
with $n_i = | \widehat{F}_i \cap V_{\gamma_i} |$ minimal.
(Recall that if $\widehat{F}_i$ is orientable then it is assumed to be essential.)
Then $n_i>0$.
Thus $F_i = \widehat{F}_i \cap M$ is a punctured surface properly embedded in $M$,
each of whose $n_i$ boundary components has slope $\gamma_i$ on $\partial_0M$.

\begin{lemma} \label{lem:incomp}
For each of six types, $F_i$ is incompressible and
boundary incompressible in $M$.
\end{lemma}

\begin{proof}
For types $\mathcal{S}$, $\mathcal{A}$ and $\mathcal{T}$,
the minimality of $n_i$ guarantees that
$F_i$ is incompressible and boundary incompressible in $M$.

For type $\mathcal{P}$, assume that $D$ is a compressing disk for $F_i$.
Since $\partial D$ is orientation-preserving on $F_i$,
$\partial D$ bounds a disk $D'$ on $\widehat{F}_i$.
Since $\mathrm{Int}\, D'$ meets $V_{\gamma_i}$,
we can create a new projective plane by replacing $D'$ with $D$,
which meets $V_{\gamma_i}$ fewer than $\widehat{F}_i$,
contradicting the minimality of $n_i$.
Next, assume that $E$ is a boundary compressing disk for $F_i$
with $\partial E = a \cup b$,
where $a$ is an essential arc in $F_i$ and $b$ lies in $\partial_0 M$.
If $a$ joins distinct components of $\partial F_i$,
then a compressing disk for $F_i$ is obtained from two parallel copies of $E$
and the disk obtained by removing a neighborhood of $b$
from the annulus in $\partial_0 M$ cobounded by those components of
$\partial F_i$ meeting $a$.
Hence both endpoints of $a$ lie in the same component,
say $\partial_1 F_i$, of $\partial F_i$.
If $n_i \geq 2$, then $b$ bounds a disk $D'$ in $\partial_0 M$
together with a subarc of $\partial_1 F_i$.
Then $E\cup D'$ gives a compressing disk for $F_i$ in $M$.
Therefore $n_i = 1$.
Then $M$ would contain a M\"{o}bius band, a contradiction.

For type $\mathcal{B}$, assume that $D$ is a compressing disk for $F_i$.
If $\partial D$, which is orientation-preserving on $F_i$,
bounds a disk $D'$ on $\widehat{F}_i$,
a new M\"{o}bius band obtained by replacing $D'$ with $D$
has fewer $n_i$, a contradiction.
Thus $\partial D$ is essential, and indeed separating on $\widehat{F}_i$.
Compressing $\widehat{F}_i$ along $D$ would give in $M(\gamma_i)$
a projective plane $Q_1$ and a disk $Q_2$ which are disjoint.
Note that the core of $V_{\gamma_i}$ must meet both $Q_1$ and $Q_2$
because $M$ is hyperbolic.
Let $\ell$ be a subarc of the core of $V_{\gamma_i}$
which connects $Q_1$ and $Q_2$, meeting them only on its endpoints.
Attaching a tube along $\ell$ to $Q_1\cup Q_2$ gives
a M\"{o}bius band in $M(\gamma_i)$ which intersects
the core of $V_{\gamma_i}$ in $n_i-2$ points, a contradiction.
Thus we have shown that $F_i$ is incompressible.
The incompressibility of $F_i$ then implies
the boundary incompressibility of $F_i$;
otherwise, $n_i = 1$ as previous,
and furthermore a booundary-compressing disk $E$
allows us to isotope the core of $V_{\gamma_i}$
into $\widehat{F}_i$ as an orientation-reversing loop.
Then $M$ would contain an essential annulus, a contradiction.

For type $\mathcal{K}$, assume that $D$ is a compressing disk for $F_i$.
As previous, $\partial D$ is
orientation-preserving and essential on $\widehat{F}_i$.
Compressing $\widehat{F}_i$ along $D$ would give in $M(\gamma_i)$
either a non-separating 2-sphere or two disjoint projective planes,
according as $\partial D$ is non-separating
or separating on $\widehat{F}_i$.
But this contradicts the irreducibility of $M(\gamma_i)$.
The incompressibility of $F_i$ then implies
the boundary incompressibility of $F_i$,
unless $M$ contains a M\"{o}bius band, as above.
\end{proof}

We use $i$ and $j$ to denote $1$ or $2$, with the convention that,
when both appear, $\{i,j\} = \{1,2\}$.

By an isotopy of $F_1$, say, we may assume that $F_1$ intersects
$F_2$ transversely. 
By Lemma \ref{lem:incomp} it can be assumed that no circle component of $F_1\cap F_2$
bounds a disk in $F_1$ or $F_2$. 
Let $G_i$ be the graph in $\widehat{F}_i$ obtained by taking
as the (fat) vertices the disks $\widehat{F}_i - \mathrm{Int}\, F_i$
and as edges the arc components of $F_i \cap F_j$ in $\widehat{F}_i$.
Thus the interior of any disk face of $G_i$ is disjoint from $F_j$.
We number the components of $\partial F_i\cap \partial_0M$ as $1,2,\cdots,n_i$
in the order in which they appear on $\partial_0 M$.
On occasion we will use $0$ instead of $n_i$ in short.
This gives a numbering of the vertices of $G_i$.
Furthermore it induces a labelling of the endpoints of
edges in $G_j$ in the usual way (see \cite{CGLS}).
Note that $G_1$ and $G_2$ have no trivial loops by Lemma \ref{lem:incomp}.
For the sake of simplicity we will say that
$F_i$ and $G_i$ are of type $X$ if $M(\gamma_i)$ is of type $X$.

Since $M$ is hyperbolic, we have the following easy lemma.

\begin{lemma} \label{lem:n}
$n_i \geq 3$ for $G_i$ of type $\mathcal{S}$, and
$n_i \geq 2$ for $G_i$ of type $\mathcal{P}$.
\end{lemma}

\begin{figure}[h]
\epsfbox{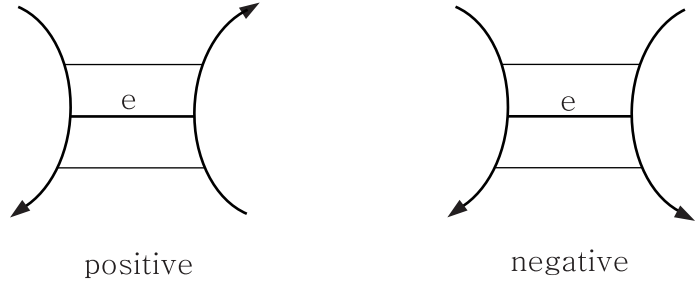}
\caption{A sign of an edge}\label{fig:parity}
\end{figure}

Although $F_i$ of type $\mathcal{P}$, $\mathcal{B}$ or
$\mathcal{K}$ is non-orientable,
we can establish a parity rule, which plays a crucial role.
In fact, this is a natural generalization of the usual parity rule \cite{CGLS}.
First, orient all components of $\partial F_i$
so that they are mutually homologous on $\partial_0 M$.
Let $e$ be an edge of $G_i$.
Since $e$ is an arc properly embedded in $F_i$,
a regular neighborhood $D$ of $e$ in $F_i$ is a disk in $F_i$.
Then $\partial D = a \cup b \cup c \cup d$,
where $a$ and $c$ are arcs in $\partial F_i$ with
induced orientations from $\partial F_i$.
If $a$ and $c$ are directed along $\partial D$,
then $e$ is called {\em positive\/}, otherwise {\em negative\/}.  See Figure \ref{fig:parity}.
Then we have the {\em parity rule}: {\it 
an edge $e$ is positive in $G_i$ if and only if
$e$ is negative in $G_j$.}

The rest of this section will be devoted to
several definitions and well known lemmas.
Let $x$ be a label of $G_i$.
An $x$-{\em edge\/} in $G_i$ is an edge with label $x$ at one endpoint,
and an $xy$-edge is an edge with label $x$ and $y$ at both endpoints.
If an $x$-edge has an endpoint at a vertex $v$ with label $x$, then it is called
an $x$-edge at $v$.
Especially, a positive $xx$-edge is called a {\em level $x$-edge\/},
which means that
the vertex $x$ of the other graph $G_j$ is incident to a negative loop.
In particular, $G_i$ does not contain a level $x$-edge unless 
$G_j$ is of type $\mathcal{P}$, $\mathcal{B}$ or $\mathcal{K}$.

An $x$-{\em cycle\/} is a cycle of positive $x$-edges of $G_i$
which can be oriented so that the head of each edge has label $x$.
A {\em Scharlemann cycle\/} is an $x$-cycle that bounds a disk face of $G_i$,
only when $n_j \geq 2$.
Each edge of a Scharlemann cycle has the same label pair $\{x,x+1\}$,
so we refer to such a Scharlemann cycle as an $(x,x+1)$-Scharlemann cycle.
The number of edges in a Scharlemann cycle
is called its {\em length}.
In particular, a Scharlemann cycle of length two is called an
{\em $S$-cycle\/} in short.

Suppose that a Scharlemann cycle $\sigma$ is immediately
surrounded by a cycle $\kappa$, that is, each edge of $\kappa$
is immediately parallel to an edge of $\sigma$.
Then $\kappa$ will be referred to as
an {\em extended Scharlemann cycle\/}, only when $n_j \geq 4$.
A {\em generalized $S$-cycle\/} is the triple
$\{e_1,e_2,e_3\}$ of mutually parallel positive edges in succession
and $e_2$ is a level edge, only when $n_j \geq 3$.

%%%%%%%%%%%%%%%%%%%%%%%%%%%%%%%%%%%%%%%%%%%%%%%%%
\begin{lemma} \label{lem:Scharsep}
Assume that $M(\gamma_j)$ is either of type $\mathcal{S}$, $\mathcal{A}$ or $\mathcal{T}$.
If $G_i$ has a Scharlemann cycle
then $\widehat{F}_j$ must be separating, and so $n_j$ is even.
Furthermore, for cases $\mathcal{A}$ and $\mathcal{T}$,
the edges of a Scharlemann cycle do no lie in a disk in $\widehat{F}_j$.
\end{lemma}

\begin{proof}
Let $E$ be a disk face bounded by a Scharlemann cycle
with a label pair, say $\{1,2\}$, in $G_i$.
Let $V_{12}$ be the $1$-handle cut from $V_{\gamma_j}$
by the vertices $1$ and $2$ of $G_j$.
(When $n_j=2$, $V_{12}$ is chosen to meet $\partial E$.)
Then tubing $\widehat{F}_j$ along $\partial V_{12}$ and
compressing along $E$ gives a new surface $R$ in $M(\gamma_j)$, homeomorphic to $\widehat{F}_j$,
that intersects $V_{\gamma_j}$ fewer times than $\widehat{F}_j$.
If the original $\widehat{F}_j$ is non-separating, then so is $R$.
Thus $R$ is essential for each case, a contradiction.

If the edges of the Scharlemann cycle lie
in a disk $D$ in an annulus or torus $\widehat{F}_j$,
then $\mbox{nhd} (D \cup V_{12} \cup E)$ is a once punctured lens space.
By the irreducibility of $M(\gamma_j)$, $M(\gamma_j)$ is a lens space, so neither of type $\mathcal{A}$ nor of type $\mathcal{T}$.
\end{proof}

\begin{lemma} \label{lem:Schar}
\begin{itemize}
\item[(1)] If $G_j$ is of type $\mathcal{S},\mathcal{A}$ or $\mathcal{T}$, 
then $G_i$ cannot have a level edge.
If $G_j$ is of type $\mathcal{P}$ or $\mathcal{B}$, then
$G_i$ has at most one label of level edges.
If $G_j$ is of type $\mathcal{K}$,
then $G_i$ has at most two labels of level edges.
\item[(2)] If $G_j$ is of type $\mathcal{P},\mathcal{B}$ or $\mathcal{K}$,
then $G_i$ cannot have a Scharlemann cycle.
If $G_j$ is of type $\mathcal{S}$ or $\mathcal{A}$,
then any two Scharlemann cycles of $G_i$ have the same label pair.
%When $G_j$ is of type $\mathcal{T}$,
%$G_i$ has at most two label pairs of Scharlemann cycles.
\item[(3)] When $G_j$ is of type $\mathcal{S}$,
$M(\gamma_j)$ contains a projective plane if $G_i$ contains an $S$-cycle.
\item[(4)] If $G_j$ is of type $\mathcal{S}$, $\mathcal{A}$ or $\mathcal{T}$,
then $G_i$ cannot have an extended $S$-cycle.
\item[(5)] If $G_j$ is of type $\mathcal{P}$, $\mathcal{B}$ or $\mathcal{K}$,
then $G_i$ cannot have a generalized $S$-cycle.
\end{itemize}
\end{lemma}

\begin{proof}
(1) Let $e$ be a negative loop based at a vertex $x$ in $G_j$.
Then $\mbox{nhd} (x\cup e)$ is a M\"{o}bius band in $\widehat{F}_j$.
Since only projective plane, M\"{o}bius band and Klein bottle
can contain at most one, one and two M\"{o}bius bands respectively,
the conclusions follow.

(2) When $G_j$ is of type $\mathcal{P}$, $\mathcal{B}$ or $\mathcal{K}$,
assume for contradiction that $G_i$ contains a Scharlemann cycle.
Then the construction in the proof of Lemma \ref{lem:Scharsep} gives
a new surface homeomorphic to $\widehat{F}_j$ in $M(\gamma_j)$, which meets
$V_{\gamma_j}$ fewer times than $\widehat{F}_j$.
This contradicts the minimality of $n_j$.

When $G_j$ is of type $\mathcal{S}$ or $\mathcal{A}$, this is \cite[Theorem 2.4]{GL1} or
\cite[Lemma 5.4(2)]{W3}, respectively.

%%%%%%%%%%%%%%
%When $G_j$ is of type $\mathcal{T}$,

(3) Let $\{e_1,e_2\}$ be an $S$-cycle in $G_i$ with label pair $\{k,k+1\}$.
Let $v_h$ be the $h$-th vertex of $G_j$ for $h=k,k+1$, and let $H$ be the part of $V_{\gamma_j}$ between
$v_k$ and $v_{k+1}$.
Let $D$ be the disk face bounded by the $S$-cycle, and
let $A$ be the M\"{o}bius band obtained by taking $H\cup D$ and shrinking $H$ radially to its core.
Since $\partial A$ is isotopic to the curve obtained from $e_1\cup e_2\cup v_k\cup v_{k+1}$ by shrinking
the vertices to points, $\partial A$ bounds a disk on $\widehat{F}_j$.
Thus $M(\gamma_j)$ contains a projective plane.

%%%%%%%%%%%%%%%
(4) These are \cite[Lemma 2.3]{W1}, \cite[Lemma 5.4(3)]{W3} and \cite[Lemma 2.10]{BZ1}, respectively.

(5) Let $\{e_1,e_2,e_3\}$ be a generalized $S$-cycle in $G_i$, where $e_2$ is a level edge with label $k$,
and $e_1, e_3$ have the same label pair $\{k-1,k+1\}$.
Let $v_h$ be the $h$-th vertex of $G_j$ for $h=k-1,k,k+1$, and
let $H$ be the part of $V_{\gamma_j}$ between the vertices $v_{k-1}$ and $v_{k+1}$ containing $v_k$.
Then $C=\mbox{nhd}\,(v_k\cup e_2)$ is a M\"{o}bius band in $\widehat{F}_j$.
Let $D$ be the disk in $F_i$ representing the parallelism of $e_1$ and $e_3$ and containing $e_2$.
Then we have a M\"{o}bius band $A$ from $H\cup D$ as before.
Note that $\mathrm{Int}\,A$ is disjoint from $\widehat{F}_j-C$.
Since $\partial A$ is isotopic to the curve obtained from $e_1\cup e_3\cup v_{k-1}\cup v_{k+1}$ by shrinking
$v_{k\pm 1}$ to points, $\partial A$ is orientation-preserving on $\widehat{F}_j$.

If $G_j$ is of type $\mathcal{P}$,
then $\partial A$ above bounds a disk $E$ in $\widehat{F}_j-C$.
Thus $A\cup E$ gives a new projective plane which meets $V_{\gamma_j}$ fewer times than $n_j$, a contradiction.

If $G_j$ is of type $\mathcal{B}$, 
then either $\partial A$ bounds a disk in $\widehat{F}_j-C$ or $\partial A$ is parallel to $\partial \widehat{F}_j$.
In the former, $M(\gamma_j)$ contains a projective plane, and hence $M(\gamma_j)$ is reducible,
contradicting the irreducibility of $M(\gamma_j)$.
In the latter, let $A_1$ be the annulus between $\partial A$ and $\partial\widehat{F}_j$.
Then $A\cup A_1$ is a M\"{o}bius band in $M(\gamma_j)$ which meets $V_{\gamma_j}$ fewer times than $n_j$.

If $G_j$ is of type $\mathcal{K}$,
then $\partial A$ bounds either a disk or a M\"{o}bius band $B$ in $\widehat{F}_j-C$.
In the former, $M(\gamma_j)$  contains a projective plane.
By the irreducibility of $M(\gamma_j)$, $M(\gamma_j)$ is the lens space $L(2,1)$.
But it is well known that $L(2,1)$ does not contain a Klein bottle \cite{BW}.
In the latter, $A\cup B$ gives a new Klein bottle in $M(\gamma_j)$ which meets $V_{\gamma_j}$
fewer times than $n_j$ again.
\end{proof}

For simplicity
we will call $x$ an $sl$-{\em label\/} of $G_i$
(or $sl$-{\em vertex\/} of $G_j$)
if $x$ is a label of either a Scharlemann cycle or a level edge
in $G_i$, according as $\widehat{F}_j$ is orientable or not.
Thus Lemma \ref{lem:Schar} implies that $G_i$ contains at most 2, 1, 2, 1 or 2 $sl$-labels
when $G_j$ is of type $\mathcal{S}$, $\mathcal{P}$, $\mathcal{A}$,
$\mathcal{B}$ or $\mathcal{K}$, respectively.
When $G_j$ is of type $\mathcal{T}$, we do not need an upper bound for the number of $sl$-labels of $G_i$, but
will show that at most three labels can be labels of $S$-cycles of $G_i$ unless $M(\gamma_j)$ contains a Klein bottle.
(See Claim \ref{claim:Tcase}.)

%Let $\mathbf{P}$ (resp. $\mathbf{N}$) denote the upper bounds for the number of
%mutually parallel positive (resp. negative) edges in $G_i$. 

\begin{lemma} \label{lem:parallel}
Let $F$ be a family of mutually parallel positive edges in $G_i$, and let $|F|$ denote the
number of edges in $F$.
\begin{itemize}
\item[(1)] If $G_j$ is of type $\mathcal{S}$ or $\mathcal{A}$,
$|F|\le\frac{n_j}{2}+1$.
Furthermore if $|F|=\frac{n_j}{2}+1$, then the first two or the last two edges of $F$ form an $S$-cycle.
\item[(2)] If $G_j$ is of type $\mathcal{P}$ or $\mathcal{B}$,
$|F|\le\frac{n_j+1}{2}$.
Furthermore if $|F|=\frac{n_j+1}{2}$,
then the first or the last edge of $F$ is a level edge.
\item[(3)] If $G_j$ is of type $\mathcal{T}$ and $n_j\ge 3$,
$|F|\le\frac{n_j}{2}+2$.
%Furthermore if $|F|>\frac{n_j}{2}$, then
%$F$ contains an $S$-cycle.
%then the first two and the last two edges of $F$ form disjoint $S$-cycles.
\item[(4)] If $G_j$ is of type $\mathcal{K}$ and $n_j\ge 2$,
$|F|\le\frac{n_j}{2}+1$.
%Furthermore if $|F|=\frac{n_j}{2}+1$ (when $n_j$ is even),
%then the first and the last edges of $F$ are level edges with distinct labels.
%And if $|F|=\frac{n_j+1}{2}$ (when $n_j$ is odd),
%then the first or the last edge of $F$ is a level edge.
\end{itemize}
\end{lemma}

\begin{proof}
(1) is \cite[Lemma 1.5(1)]{W2} and \cite[Lemma 2.5(2)]{GW}.
(We remark that when $G_j$ is of type $\mathcal{A}$ and $n_j=1$, $|F|\le 1$, because
a pair of edges cannot be parallel in both graphs \cite[Lemma 2.1]{G1}.)
(3) is \cite[Lemma 1.4]{W2}.
For (2), if $|F|>\frac{n_j+1}{2}$ then $F$ would contain either an $S$-cycle, a generalized $S$-cycle or
two level edges with distinct labels, except the case where $G_j$ is of type $\mathcal{B}$ and $n_j=1$.
In the exceptional case, a pair of edges is parallel in both $G_i$ and $G_j$,
a contradiction.
(4) is similar to (2).
\end{proof}

\begin{lemma}\label{lem:parallel-negative}
Let $F$ be a family of mutually parallel negative edges in $G_i$.
\begin{itemize}
\item[(1)] When $G_j$ is of type $\mathcal{S}$ or $\mathcal{P}$,
           $|F|\le n_j-1$.
\item[(2)] When $G_j$ is of type $\mathcal{A}$,
           $\mathcal{B}$ or $\mathcal{K}$,
           $|F|\le n_j$.
\end{itemize}
\end{lemma}

\begin{proof}
(1) See \cite[Lemma 2.3]{G1}.  (This is essentially proved in \cite[Section 5]{GLi}.)
(2) The argument of the proof of \cite[Lemma 2.5(3)]{GW} works well.  (Also, see the proof of \cite[Lemma 4.2]{G1}.)
\end{proof}

%%%%%%%%%%%%%%%%%%%%%%%%%%%%%%%%%%%%%%%%
\section{$x$-faces and non-orientable surfaces} \label{sec:Pxface}

In this section we assume that $G_j$ is of type $\mathcal{P}$ or $\mathcal{B}$
(so we may assume that $1$ is the only possible $sl$-vertex of $G_j$ by Lemma \ref{lem:Schar}),
and thus all level edges of $G_i$ are level $1$-edges.

A disk face of the subgraph of $G_i$ consisting of all the vertices
and positive $x$-edges of $G_i$ is called an {\em $x$-face\/}.
Remark that the boundary of an $x$-face $D$ may be not a circle, that is,
$\partial D$ may contain a double edge, and more than two edges of $\partial D$
may be incident to a vertex on $\partial D$
(see Figure 2.1 in \cite{HM}).
A cycle in $G_i$ is a {\em two-cornered cycle\/}
if it is the boundary of a disk face containing only $01$-corners,
$12$-corners and positive edges, only when $n_j \geq 3$.
Recall that $0$ denotes $n_j$.
A two-cornered cycle must contain both kinds of corners,
because $G_i$ cannot contain a Scharlemann cycle by Lemma \ref{lem:Schar}.
Also, it contains at least one $02$-edge.
For convenience, when the labels appear in anticlockwise order
around the boundary of a vertex $v$ of $G_i$,
given three distinct labels $x_1,x_2,x_3$,
we say $x_1 < x_2 < x_3$ if $x_1,x_2,x_3$ appear
in anticlockwise order on some interval in $\partial v$ containing $n_j$ edge endpoints.
Thus three expressions $x_1<x_2<x_3$, $x_2<x_3<x_1$ and $x_3<x_1<x_2$ are equivalent.

\begin{proposition} \label{prop:Ppair}
Suppose that $n_j \geq 3$.
An $x$-face, $x \ne 1$, in $G_i$ contains a pair of two-cornered
cycles sharing a level $1$-edge.
\end{proposition}

\begin{proof}
Let $\Gamma_D$ be the subgraph of $G_i$ in an $x$-face $D$.
There is a possibility that $\partial D$ is not a circle as mentioned before.
Since we will find a pair of two-cornered cycles within $D$,
we can cut formally the graph $G_i \cap D$ along double edges of $\partial D$
and at vertices to which more than two edges of $\partial D$ are incident to
so that $\partial D$ is deformed into a circle.
(See also Figure 5.1 in \cite{HM}.)
Thus we may assume that $\partial D$ is a circle.
%$\Gamma_D$ has no vertex in the interior of $D$.
We may assume that the labels appear in anticlockwise order
around the boundary of each vertex.

Suppose that $D$ has a diagonal edge $d$ with distinct labels
$\{a,b\}$, which must differ from $x$ because $D$ is an $x$-face,
as in Figure \ref{fig:split}(a).
Assume without loss of generality that $a<x<b$.
Formally construct a new $x$-face $D'$ as follows.
Keep all corners and edges of $\Gamma_D$ to the right of $d$
(when $d$ is directed from $a$ to $b$),
discard all corners and edges to the left of $d$,
and then insert additional edges to the left of $d$,
and parallel to $d$, until you first reach label $x$
at one or both ends of this parallel family of edges,
as in Figure \ref{fig:split}(b).
In particular, these additional edges contain no edges of two-cornered cycles
or Scharlemann cycles of the graph on the new $x$-face $D'$.

\begin{figure}[h]
\epsfbox{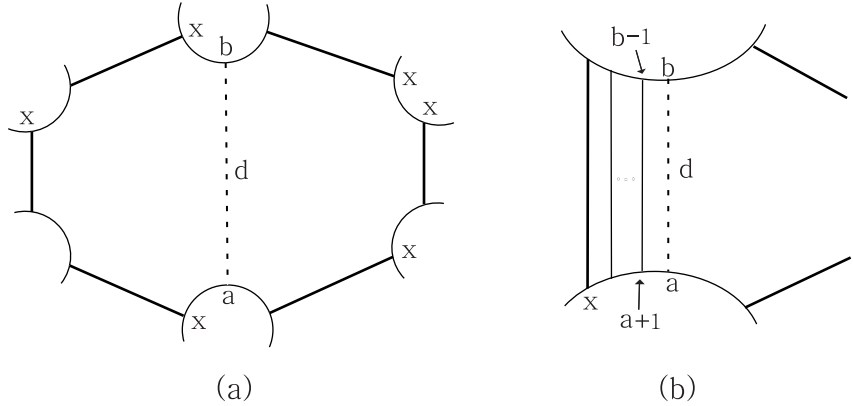}
\caption{Split along a diagonal edge}\label{fig:split}
\end{figure}

Repeat the above process for every diagonal edge which is not a level $1$-edge,
then get a new $x$-face $E$ and a graph $\Gamma_E$ in $E$.
All diagonal edges of $\Gamma_E$ are level $1$-edges,
and all (and only) boundary edges are $x$-edges,
where the label $x$ possibly appear on both ends,
say level $x$-edges.

\begin{claim} \label{claim:level}
$\Gamma_E$ contains a level $1$-edge.
\end{claim}

\begin{proof}[Proof of Claim \ref{claim:level}]
Assume that $\Gamma_E$ contains no level $1$-edges,
and so no diagonal edges.
We first show that if for some vertex $v$ of $\Gamma_E$
two boundary edges are incident to $v$ with label $x$,
then these should be level $x$-edges.
For, $n_j + 1$ edges are incident to $v$ in $\Gamma_E$.
If $n_j$ is even, more than $\frac{n_j}{2}$ mutually parallel edges
are incident to $v$, and so one of these edges should
be a level edge different from a level $x$-edge
(recall that $\Gamma_E$ cannot contain a Scharlemann cycle),
a contradiction.
If $n_j$ is odd, two families of $\frac{n_j+1}{2}$ mutually parallel
edges are incident to $v$ by Lemma \ref{lem:parallel}(2), and the boundary edge of each family
should be a level $x$-edge by the same reason above.

Consider the cycle $\sigma$ consisting of boundary $x$-edges of $\Gamma_E$.
Assume that $\sigma$ has an $x$-edge which is not a level $x$-edge.
So the only one end has label $x$ at $v_1$, say.
By the fact we just proved,
another $x$-edge incident to $v_1$ does not have label $x$ at $v_1$.
Thus this edge has label $x$ at the other end $v_2$, say.
After repeating this process,
we are led to show that $\sigma$ is a great $x$-cycle
in the terminology of \cite{CGLS}.
By the same argument in the proof of Lemma 2.6.2 of \cite{CGLS},
$\Gamma_E$ contains a Scharlemann cycle,
since $\Gamma_E$ does not contain level edges in its interior,
a contradiction.

If all edges of $\sigma$ are level $x$-edges, then we have a great $x+1$-cycle just inside $\sigma$.
Thus we still find a Scharlemann cycle, a contradiction.
(See also \cite[Lemma 5.2]{HM}.)
\end{proof}

Let $e$ be a level $1$-edge.
So it does not belong to $\partial \Gamma_E$.
Let $E_1$ and $E_2$ be the faces of $\Gamma_E$ adjacent to $e$.
Assume for contradiction that $\partial E_1$ (or $\partial E_2$)
is not a two-cornered cycle.
Note that $\partial E_1$ may contain many level $1$-edges.
Let $\{a_k,a_k+1\}$, $k=1, \cdots ,n$, be the consecutive label pairs
of the corners between successive level $1$-edges on $\partial E_1$, which appear in order around $\partial E_1$,
when one runs clockwise around $\partial E_1$
starting at one end of $e$, as in Figure \ref{fig:cornered}(a).
Then some $a_k$ is neither $0$ nor $1$.
Since $a_1 = 1$ and $a_n = 0$,
there are indices $l$ and $m$ so that
$a_k = 0$ or $1$ when $1 \leq k < l$ or $k=m$,
and $a_k \neq 0,1$ when $l \leq k <m$.
Consider the edges of the parallelism class containing
each $\{a_{k-1} +1,a_k \}$-edge, $l \leq k \leq m$.
Note that among these edges a level $x$-edge on the boundary
of $\Gamma_E$ is the only possible level edge,
and there is no $x$-edge except on the boundary.
See Figure \ref{fig:cornered}(b).
Then we have $x \leq a_k < a_{k-1}+1 \leq x$, and so $x \leq a_k \leq a_{k-1} < x$.
Finally we have
$x \leq a_m \leq a_{m-1} \leq \cdots \leq a_l \leq a_{l-1} < x$.
This is impossible because $a_{l-1}, a_m = 0$ or $1$.
Thus we have shown that any face adjacent to a level $1$-edge is two-cornered.
%
%Choose an outermost level $1$-edge in $\Gamma_E$.
%Then two faces adjecent to it are a desired pair of two-cornered cycles.
\end{proof}

In fact, we can see from the argument above that 
$a_1=\dots=a_i=1$ and $a_{i+1}=\dots=a_n=0$ for some $i$.
That is, when we go along a two-cornered cycle from a level $1$-edge in some direction, 
$12$-corners appear successively, $23$-corners appear after them, and
we reach a level $1$-edge (possibly, a different one from the start).

\begin{figure}[h]
\epsfbox{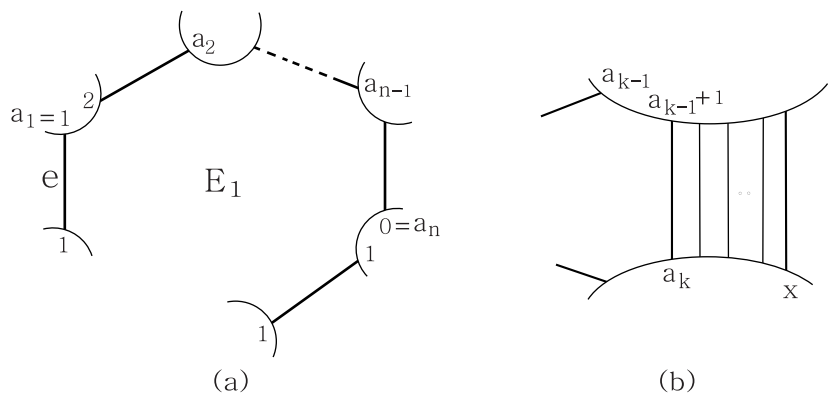}
\caption{Finding two-cornered cycles}\label{fig:cornered}
\end{figure}

Recall that $1$ is the only possible $sl$-label of $G_i$.
The next theorem essentially follows from the argument of \cite[Section 7]{H}.

\begin{theorem} \label{thm:Pxface}
If $G_j$ is of type $\mathcal{P}$,
then $G_i$ cannot contain an $x$-face for a non-$sl$-label $x$ of $G_i$.
\end{theorem}

\begin{proof}
Assume that $G_i$ contains such an $x$-face.
If $n_j = 2$, then
a boundary of each face in the $x$-face is a Scharlemann cycle
or contains both a level $1$-edge and a level $2$-edge,
contradicting Lemma \ref{lem:Schar}(1) and (2).
Thus we may assume that $n_j \geq 3$.
Applying Proposition \ref{prop:Ppair},
$G_i$ contains a pair of two-cornered cycles sharing a level $1$-edge $e$.

Let $\widehat{S}$ be the sphere which is the boundary of a regular neighborhood
of $\widehat{F}_j$ in $M(\gamma_j)$, and let $S=\widehat{S}\cap M$.
The arc components of $F_i \cap S$
give rise to a pair of labelled graphs ($G_i^S, G_S$)
as usual, where $G_S$ is a double cover of $G_j$.
Then $G_i^S$ contains an $S$-cycle $\sigma$ corresponding to the level $1$-edge $e$.
We may assume that $\sigma$ has label pair $\{1,2\}$.
Also, two cycles adjacent to $\sigma$ in $G_i^S$ contain
only two types of corners, $01$-corners and $23$-corners.
Let $E_1$ and $E_2$ be the faces of $G_i^S$ bounded by
these cycles respectively.
The edges of $\sigma$ separate $\widehat{S}$ into two disks,
and two vertices $0$ and $3$ of $G_S$ are contained
in the same disk component, say $D_0$,
because of the existence of $03$-edges.
Let $D_1$ be an expansion of $D_0$
containing the fat vertices $1$ and $2$ in $\widehat{S}$,
and regard it as properly embedded in
$X = \mbox{nhd}(D_1 \cup V_{01} \cup V_{23} \cup E_1 \cup E_2)$
where $V_{01}$ and $V_{23}$ are defined
as in the proof of Lemma \ref{lem:Scharsep}.
Since $\partial E_1$ and $\partial E_2$ are non-separating on the boundary of
the genus two handlebody $\mbox{nhd}(D_1 \cup V_{01} \cup V_{23})$,
$\partial X$ is either a $2$-sphere or
the disjoint union of a $2$-sphere and a torus.
The latter happens only when $\partial E_1$ and $\partial E_2$ are parallel on
the boundary of $\mbox{nhd}(D_1 \cup V_{01} \cup V_{23})$.
In this case, $\partial D_1$ lies in the parallelism annulus.
Thus $\partial D_1$ cuts the $2$-sphere component of $\partial X$
into two disk parts, one of which is parallel to $D_1$, in any case.
Let $D$ be the other disk part and
then extend $\partial D$ down to $\widehat {F}_j$
along an annulus $A = \partial D \times I$,
using the $I$-bundle structure of $\mbox{nhd}(\widehat {F}_j)$.
Note that the boundary component of $A$
which lies on $\widehat {F}_j$ bounds a M\"{o}bius band $B$
intersecting the core of $V_{\gamma_j}$ in a single point.
In particular,
$D \cup A \cup B$ is a projective plane in $M(\gamma_j)$.
Since $|D\cap V_{\gamma_j}|\leq |D_0\cap V_{\gamma_j}|-2=n_j -3$,
$|A \cap V_{\gamma_j}| = 0$ and $|B \cap V_{\gamma_j}| = 1$,
this projective plane intersects the core of $V_{\gamma_j}$
in less than $n_j$ points.
This contradicts the minimality of $n_j$.
\end{proof}

\begin{theorem} \label{thm:Bxface}
If $G_j$ is of type $\mathcal{B}$ with $n_j \geq 2$,
then $G_i$ cannot contain an $x$-face for a non-$sl$-label $x$ of $G_i$.
\end{theorem}

\begin{proof}
The same argument in the proof of Theorem \ref{thm:Pxface}
applies here.
Assume that $G_i$ contains an $x$-face for a non-$sl$-label $x$.
Applying Proposition \ref{prop:Ppair},
$G_i$ contains a pair of two-cornered cycles sharing a level $1$-edge $e$.
In particular, we can assume that one of the pair contains a single level $1$-edge.
This is guaranteed by choosing the level $1$-edge $e$ to be outermost in the $x$-face.
In the same notation as before, 
however, the edges of $\sigma$ separate $\widehat{S}$, which is an annulus,
into two annuli, and so $D_1$ is an annulus.
Also, we may assume that $\partial E_1$ contains only one $12$-edge.

Let $W=\mbox{nhd}(D_1\cup V_{01}\cup V_{23})$ and let
$X=W\cup \mbox{nhd}(E_1 \cup E_2)$ as before.
We may assume that $D_1$ is properly embedded in $W$ and $X$.
Thus $W$ consists of a solid torus $\mbox{nhd}(D_1)$ and the two $1$-handles attached there.
Now $\pi_1(W)=\langle x,y,t \rangle$, where taking as base \lq\lq point\rq\rq\ a
thick ball negihborhood of a meridian disk of $\mbox{nhd}(D_1)$, which contains
the $12$-edge of $\partial E_1$ and the four attaching disks
of two $1$-handles $V_{01}$ and $V_{23}$,
$x$ is represented by a core of $V_{23}$ going from vertex $2$ to vertex $3$,
$y$ is represented by a core of $V_{01}$ going from vertex $0$ to vertex $1$, and
$t$ is represented by the $12$-edge of $\sigma$ on $\partial E_2$ oriented from vertex $1$ to vertex $2$.
In $\pi_1(W)$, $[\partial E_2]$ (with a clockwise orientation) contains the sequence $ytx$,
but $[\partial E_1]$ has no such sequence by the fact that
$\partial E_1$ has only one level $12$-edge.
Hence $\partial E_1$ and $\partial E_2$ are not parallel on $\partial W$.
Clearly, both are non-separating on $\partial W$.
Thus the boundary of $W\cup \mbox{nhd}(E_1)$ is a genus two surface, but
$\partial X$ is either a torus or the disjoint union of two tori, according to
whether $\partial E_2$ is non-separating on the genus two surface or not.
In the former, one component of $\partial D_1$ is essential on the torus.
If not, $\widehat{F}_j$ is compressible, and hence
$M(\gamma_j)$ contains a projective plane, which contradicts the irreducibility of $M(\gamma_j)$.
Therefore the frontier of $X$ is an annulus.
In the latter, the frontier of $X$ is the disjoint union of an annulus and a torus.
Thus in either case, we have an annulus component among the frontier of $X$.
In particular, $\partial D_1$ divides it into
two annuli, one of which is parallel to $D_1$.
The rest of the proof is exactly the same as previous, and
we would have a new M\"{o}bius band having fewer intersection with $V_{\gamma_j}$.
\end{proof}

%Remark that a generalized $S$-cycle is one of the simplest form of $x$-faces.
%Thus two theorems above guarantee that
%if $G_j$ is of type $\mathcal{P}$ or $\mathcal{B}$,
%then $G_i$ cannot contain a generalized $S$-cycle.

%%%%%%%%%%%%%%%%%%%%%%%%%%%%%%%%%%%%%%%%%%%%%%%%%%%%%%%%%%%%%%%%%%%%%%%%%%%%%%%%%%%%%%%
\section{$x$-faces and orientable surfaces} \label{sec:Sxface}

In this section we assume that $G_j$ is of type $\mathcal{S}$ or $\mathcal{A}$.
By Lemma \ref{lem:Schar}(2) we may assume that
$G_i$ contains only $12$-Scharlemann cycles, if they exist.

An $x$-face in $G_i$ is defined as previous.
A cycle in $G_i$ is a {\em two-cornered cycle\/},
slightly different from previous,
if it is the boundary of a disk face containing only $01$-corners,
$23$-corners and positive edges, and additionally
it contains at least one edge of a $12$-Scharlemann cycle.
A two-cornered cycle must contain both types of corners and a $03$-edge.
A {\em cluster\/} $C$ is a connected subgraph of $G_i$ satisfying that
\begin{itemize}
\item[(i)] $C$ consists of $12$-Scharlemann cycles and two-cornered cycles,
\item[(ii)] every $12$-edge of $C$ belongs to both a Scharlemann cycle
and a two-cornered cycle, and
\item[(iii)] $C$ contains no cut vertex.
\end{itemize}

\begin{figure}[h]
\epsfbox{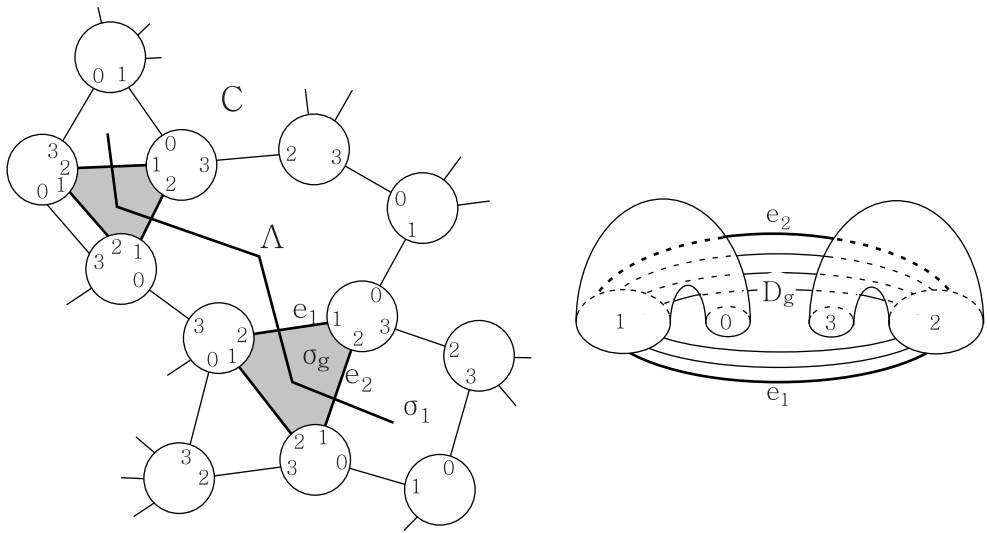}
\caption{A cluster and a seemly pair}\label{fig:cluster}
\end{figure}

See Figure \ref{fig:cluster}.

\begin{proposition} \label{prop:cluster}
Suppose that $n_j \geq 3$.
An $x$-face, $x \neq 1,2$, in $G_i$ contains a cluster $C$.
\end{proposition}

\begin{proof}
Let $\Gamma_D$ be the subgraph of $G_i$ in an $x$-face $D$.
As in the previous,
we may assume that $\partial D$ is a circle.
%and $\Gamma_D$ has no vertex in the interior of $D$.
Also assume that the labels appear anticlockwisely.

For all diagonal edges of $D$ which are not of $12$-Scharlemann cycles
(also these are neither $x$-edges nor level edges),
apply the same argument in the proof of Proposition \ref{prop:Ppair}.
Then we get a new $x$-face $E$ and a graph $\Gamma_E$ in $E$
so that all diagonal edges are of $12$-Scharlemann cycles
and all (and only) boundary edges are $x$-edges.
Furthermore the additional edges contain no edges of
Scharlemann cycles or two-cornered cycles of $\Gamma_E$.
Remark that a level $x$-edge can appear
on the boundary of the graph.

\begin{claim} \label{claim:Schar}
$\Gamma_E$ contains a $12$-Scharlemann cycle,
so does $\Gamma_D$.
\end{claim}

\begin{proof}[Proof of Claim \ref{claim:Schar}]
It is clear from the proof of Claim \ref{claim:level}.
\end{proof}

This means that $\widehat{F}_j$ must be separating
and $n_j$ is even by Lemma \ref{lem:Scharsep}.
The parity rule guarantees that each edge of $\Gamma_E$
connects vertices with one label even and the other label odd,
and so there are no level $x$-edges.

Any $12$-edge of a Scharlemann cycle
does not belong to $\partial \Gamma_E$.
Consider the face $E_1$ of $\Gamma_E$ which is adjacent to the $12$-edge
and whose boundary is not the Scharlemann cycle.
It is possible that $E_1$ contains more than one $12$-edges
of Scharlemann cycles.
Again, let $\{a_k,a_k+1\}$, $k=1, \cdots ,n$,
be the consecutive label pairs of the corners
between two consecutive $12$-edges of Scharlemann cycles
when one runs clockwise around $\partial E_1$.
Note that $a_1 = 2$ and $a_n = 0$.

Assume for contradiction that
$\partial E_1$ is not a two-cornered cycle.
Since some $a_k$ then is neither $0$ nor $2$,
there are indices $l$ and $m$ so that
$a_k = 0$ or $2$ when $1 \leq k < l$ or $k=m$,
and $a_k \neq 0,2$ when $l \leq k <m$.

Consider the edges of the parallelism class containing
each $\{a_{k-1} +1,a_k \}$-edge for $l \leq k \leq m$.
Since there are neither Scharlemann cycles
nor level edges among these edges,
one finds that $x \leq a_k < a_{k-1}+1 \leq x$, and hence
$x \leq a_k \leq a_{k-1} < x$.
And so
$x \leq a_m \leq a_{m-1} \leq \cdots \leq a_l \leq a_{l-1} < x$.
This is impossible because $a_{l-1}, a_m = 0$ or $2$
and all $a_k$'s are even by the parity rule.
Hence $\partial E_1$ is a two-cornered cycle.

Thus we have shown that 
any face next to a Scharlemann cycle in $\Gamma_E$ is two-cornered.
Let $C$ be the union of all the Scharlemann cycles and
all the two-cornered cycles adjacent to each 12-edges
of the Scharlemann cycles.
If necessary, choose a block of $C$. Then it is a desired cluster
in $\Gamma_E$ and so in $\Gamma_D$.
\end{proof}

Furthermore, we can see from the argument that
$a_1=\dots=a_i=2$ and $a_{i+1}=\dots=a_n=0$ for some $i$.

Let $R$ be the twice-punctured sphere obtained from
$\widehat{F}_j$ by deleting two fat vertices $1$ and $2$
(if $G_j$ is of type $\mathcal{A}$, then use $\widehat{F}_j$ after capping
off two boundary circles by disks).
The family of all $12$-edges of a Scharlemann cycle in the cluster $C$
separates $R$ into disks,
and one of those disks contains both vertices $0$ and $3$ of $G_j$,
because of the existence of $03$-edges in $C$.
The two $12$-edges bounding such a disk are
called {\em good edges\/} of $C$.
Thus each Scharlemann cycle in $C$ has exactly two good edges.

Let $\Lambda$ be the maximal dual graph of $C$
whose vertices are dual to Scharlemann cycles
and two-cornered cycles containing good edges,
and edges are dual to good edges of $C$
as depicted in Figure \ref{fig:cluster}.
Thus in $\Lambda$, a vertex dual to a Scharlemann cycle has valency $2$,
and a vertex dual to a two-cornered cycle has valency
the number of good edges of the two-cornered cycle.
Furthermore $\Lambda$ is a forest according to the construction of $C$.
That is, each component of $\Lambda$ is a tree.
Let $\sigma$ be a two-cornered cycle with a good edge in $C$.
Then each $12$-edge, which is not a good one, in $\sigma$
contributes to the number of components of $\Lambda$ by adding $1$.
Consequently there is a component $\Lambda_g$ of $\Lambda$ so that
all $12$-edges of its dual two-cornered cycles are good.

Hereafter, we consider the subgraph $C_g$ of $C$, dual to $\Lambda_g$.
Say, $C_g$ contains $n$ Scharlemann cycles,
and so $2n$ good edges and $n+1$ two-cornered cycles.
A two-cornered cycle dual to an end vertex of the tree $\Lambda_g$
has only one good edge.
Choose one $e_1$ of the nearest edges to vertex $0$ (or $3$) among them,
that is, there are no such good edges
between $e_1$ and vertex $0$ in $R$.
Let $\sigma_g$ and $\sigma_1$ be the Scharlemann cycle and
two-cornered cycle adjacent to $e_1$ respectively.
Note that $\sigma_1$ has only one $12$-edge.
Then $\sigma_g$ has another good edge $e_2$,
and $e_1$ and $e_2$ bound a disk $D_g$
containing $0$ and $3$ in $R$.
Note that the boundaries of the faces bounded by
Scharlemann cycles are parallel on a torus obtained
from $R$ by attaching an annulus $\partial V_{12}$.
Thus exactly $n-1$ out of $2n$ good edges are not contained in $D_g$.
Therefore we have another two-cornered cycle $\sigma_2$
all of whose $12$-edges lie in $D_g$.
Consequently all edges (consisting of $01$-edges, $12$-edges,
$23$-edges and $03$-edges) of $\sigma_1$ and $\sigma_2$ lie in $D_g$.
Furthermore if $\sigma_2$ has only one good $12$-edge,
then the two good edges of $\sigma_1$ and $\sigma_2$ lie on different
sides of the vertices $0$ and $3$ in $D_g$ by the choice of $\sigma_1$.
Such $\sigma_1, \sigma_2$ are called a {\em seemly pair\/}.
Then we can say:

\begin{proposition} \label{prop:seemlypair}
There is a seemly pair of two-cornered cycles in $C$.
\end{proposition}

From now we apply the argument in \cite[Section 6]{H} to get
the following two theorems.
Recall that $1,2$ are the only possible $sl$-labels of $G_i$

\begin{theorem} \label{thm:Sxface}
If $G_j$ is of type $\mathcal{S}$,
then $G_i$ cannot contain an $x$-face for a non-$sl$-label $x$ of $G_i$.
\end{theorem}

\begin{proof}
Suppose that $G_i$ contains such an $x$-face.
Note that $n_j \geq 3$ by Lemma \ref{lem:n}.
We continue the preceding argument.

Let $E$ be the disk bounded by $\sigma_g$.
Then $L=\mathrm{nhd}((\widehat{F}_j-\mathrm{Int}\,D_g)\cup V_{12}\cup E)$ is a punctured lens space in $M(\gamma_j)$.
Since $\partial L$ is a reducing sphere in $M(\gamma_j)$,
$|\partial L\cap V_{\gamma_j}|=2|(\widehat{F}_j-\mathrm{Int}\,D_g)\cap V_{\gamma_j}|-2\ge n_j$
by the minimality of $n_j$.
Thus $|\mathrm{Int}\,D_g\cap V_{\gamma_j}|\le \frac{n_j}{2}-1$.
Let $X_D' = \mbox{nhd} (D_g \cup V_{01} \cup V_{23})$ and
$X_F' = \mbox{nhd} (\widehat{F}_j \cup V_{01}\cup V_{23})$.
Then $X_D'$ is a genus two handlebody and
$X_F'$ is a once-punctured genus two handlebody.
The genus two torus component of $\partial X_F'$
is referred to as the outer boundary of $X_F'$.
Let $E_i$ be the face bounded
by the two-cornered cycle $\sigma_i$ for $i=1,2$.
The point is that all edges of $\sigma_i$ are contained in $D_g$.
Let $X_D = X_D' \cup \mbox{nhd}(E_1\cup E_2)$
and $X_F = X_F' \cup \mbox{nhd}(E_1\cup E_2)$
as in Figure \ref{fig:XdXf}.
Since $\partial E_1$ and $\partial E_2$ are non-separating on $\partial X_D'$ and $\partial X_F'$,
both $\partial X_D$ and the outer components of $\partial X_F$ are
either a $2$-sphere or the disjoint union
of a $2$-sphere and a torus, simultaneously.
Note that the latter case occurs
only when $\partial E_1$ and $\partial E_2$ are parallel on the boundary of $X_D'$ or $X_F'$.

\begin{figure}[h]
\epsfbox{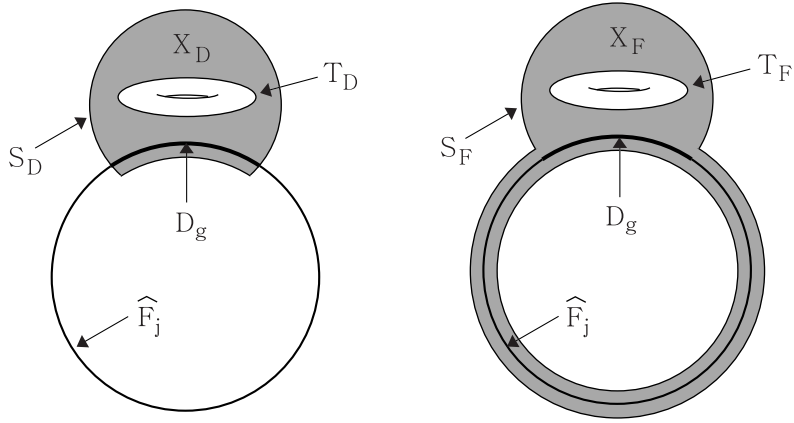}
\caption{$X_D$ and $X_F$}\label{fig:XdXf}
\end{figure}

First, assume that $\partial X_D$ is a $2$-sphere $S_D$,
and the outer component of $\partial X_F$ is a $2$-sphere $S_F$.
If $X_D$ is not a $3$-ball, then $S_D$ is a reducing sphere.
But $|S_D \cap V_{\gamma_j}|=2|\mathrm{Int}\,D_g \cap V_{\gamma_j}| \leq n_j-2$,
contradicting the minimality of $n_j$.
Thus it should be a $3$-ball,
and so $X_F$ is homeomorphic to $S^2\times I$.
Thus $S_F$ is isotopic to $\widehat{F}_j$,
contradicting the minimality of $n_j$ again.

Next, assume that $\partial X_D$ is
the disjoint union of a $2$-sphere and a torus, $S_D \cup T_D$ and
the outer components of $\partial X_F$ are also
the disjoint union of a $2$-sphere and a torus, $S_F \cup T_F$.
Recall that this case occurs whenever $\partial E_1$ and $\partial E_2$
cobound an annulus in $\partial X_D'$ and $\partial X_F'$.
This is possible only when $\sigma_2$ corresponds to
an end vertex of $\Lambda_g$ with the same number of
$01$-corners and $23$-corners as those of $\sigma_1$.
Let $D'$ be the intersection of $\partial X_D$ and
the inner sphere component of $\partial X_F$.
By the choice of the seemly pair,
this annulus in $\partial X_D'$ contains $D'$, so does $S_D$.
Similarly $S_F$ contains a pushoff of $R- D_g$.

If $S_D$ is non-separating in $M(\gamma_j)$,
then it is a reducing sphere with
$|S_D \cap V_{\gamma_j}| \leq n_j-2$, contradicting.
Thus $S_D$ is separating in $M(\gamma_j)$.
Let $X_D''$ be the manifold bounded by $S_D$
containing $T_D$ in $M(\gamma_j)$.
If $X_D''$ is not a $3$-ball,
then $S_D$ is a reducing sphere with less intersection with $V_{\gamma_j}$.
Thus it should be a $3$-ball,
and so $S_F$ is isotopic to $\widehat{F}_j$, contradicting again.
This completes the proof.
\end{proof}

%%%%%%%%%%%%%%%%%%%%%%%%%%%%%%%%%%%%%%%%%%%%%%%%%%%%%%%%%%%%%%%%%%%%%%%%%%%%%%
\begin{theorem} \label{thm:Axface}
If $G_j$ is of type $\mathcal{A}$ with $n_j \geq 3$,
then $G_i$ cannot contain an $x$-face for a non-$sl$-label $x$ of $G_i$.
\end{theorem}

\begin{proof}
Suppose that $G_i$ contains such an $x$-face.
Again we continue the argument stated before Proposition \ref{prop:seemlypair}.
The proof is divided into three cases.
Recall that $R$ is a twice-punctured sphere obtained from $\widehat{F}_j$ by deleting
two vertices $1$ and $2$ and then capping off $\partial \widehat{F}_j$ by two disks,
and that $D_g$ lies on $R$.

%%%%%%%%%%%%%
First, suppose that $D_g$ contains $\partial \widehat{F}_j$.
Then the edges of the Scharlemann cycle $\sigma_g$ lie
in a disk in $\widehat{F}_j$, contradicting Lemma \ref{lem:Scharsep}.

%%%%%%%%%%%%%%
Second, suppose that $D_g$ does not contain any component
of $\partial \widehat{F}_j$.
Then $D_g$ is contained in $\widehat{F}_j$.
The proof is similar to that of the preceding theorem.
Recall that $E_i$ is the face bounded
by the two-cornered cycle $\sigma_i$ for $i=1,2$.
Let $X_D = \mbox{nhd} (D_g \cup V_{01} \cup V_{23} \cup E_1 \cup E_2)$
and $X_F = \mbox{nhd} (\widehat{F}_j \cup V_{01} \cup V_{23}\cup E_1 \cup E_2)$.
Then $\partial X_D$ is either a $2$-sphere or the disjoint union
of a $2$-sphere and a torus.
Also, the frontier of $X_F$ has a copy of $\widehat{F}_j$, and
the remaining part is referred to as the outer components.
It can be seen that
the outer components consist of either a single annulus or
the disjoint union of an annulus and a torus by the choice of the seemly pair.
($\partial X_D$ is disconnected if and only if the outer components of the frontier of $X_F$ are disconnected.)
In either case, we have an annular component $A_F$.
Then $A_F$ is obtained from $\widehat{F}_j$
by replacing $D_g$ with a proper disk
which is a part of the $2$-sphere component of $\partial X_D$.
Since $M(\gamma_j)$ is irreducible by the assumption in Section \ref{sec:intro},
the $2$-sphere component of $\partial X_D$ bounds a $3$-ball.
Hence $A_F$ is isotopic to $\widehat{F}_j$.
This contradicts the minimality of $n_j$.

%%%%%%%%%%%%%%%%%%%%
Finally suppose that $D_g$ contains exactly one component
of $\partial \widehat{F}_j$.
If there is a disk in $D_g$ which contains all edges of $\sigma_1$ and $\sigma_2$ and which does not contain
the component of $\partial\widehat{F}_j$, then the case reduces to the previous one.
Otherwise, let $A_1$ be the annulus $D_g \cap \widehat{F}_j$.
Consider $Y=\mathrm{nhd}((\widehat{F}_j-\mathrm{Int}\,A_1)\cup V_{12}\cup E)$,
where $E$ is the face bounded by $\sigma_g$.
Then the frontier $Q$ of $Y$ in $M(\gamma_j)$ is an essential annulus by Claim on page 430 of \cite{W3}.
Thus $|Q\cap V_{\gamma_j}|=2|(\widehat{F}_j-\mathrm{Int}\,A_1)\cap V_{\gamma_j}|-2\ge n_j$
by the minimality of $n_j$.
Hence $|\mbox{Int}\,A_1 \cap V_{\gamma_j}| \leq \frac{n_j}{2} - 1$.
Let $X_A = \mbox{nhd} (A_1 \cup V_{01} \cup V_{23} \cup E_1 \cup E_2)$
and $X_F = \mbox{nhd} (\widehat{F}_j \cup V_{01} \cup V_{23}
\cup E_1 \cup E_2)$ again.
Then the frontier of $X_A$ contains an annulus component $A_A$ by
the same argument as in the proof of Theorem \ref{thm:Bxface}.
(By the choice of the seemly pair, $\partial E_1$ and $\partial E_2$ are not parallel on the boundary of
$\mbox{nhd} (A_1 \cup V_{01} \cup V_{23})$ and $\mbox{nhd} (\widehat{F}_j \cup V_{01} \cup V_{23})$.)
Also the outer components of the frontier of $X_F$ contain an annulus $A_F$.
If $A_A$ is essential in $M(\gamma_j)$,
then $|A_A \cap V_{\gamma_j}| \le 2|\mbox{Int}\,A_1 \cap V_{\gamma_j}| \leq n_j-2$,
contradicting the minimality of $n_j$.
Remark that $A_A$ is incompressible, because
the central curve of $A_A$ is isotopic to the central curve of
$\widehat{F}_j$.
Thus $A_A$ should be boundary parallel.
That is, there is a solid torus $U$ such that $A_A$ is a longitudinal annulus on $\partial U$ and
$\partial U-A_A\subset \partial M(\gamma_j)$.
Then either $\widehat{F}_j$ is boundary parallel or
one can isotope $\widehat{F}_j$ through $U$ to $A_F$, which intersects $V_{\gamma_j}$ fewer times than $\widehat{F}_j$.
\end{proof}

Remark that an extended Scharlemann cycle is one of the simplest form of $x$-faces.
Thus two theorems above guarantee that
if $G_j$ is of type $\mathcal{S}$ or $\mathcal{A}$,
then $G_i$ cannot contain an extended Scharlemann cycle.
We may emphasize that only an extended $S$-cycle was used
in the literatures such as \cite{GW,W1} etc.

%%%%%%%%%%%%%%%%%%%%%%%%%%%%%%%%%%%%%%%%%%%%%%%%%%%%%%%%%
\section{Extremal block} \label{sec:extremal}

The {\em reduced graph\/} $\overline{G}_i$ of $G_i$
is defined to be the graph obtained from $G_i$ by
amalgamating each family of parallel edges into a single edge.
Let $G_i^+$ denote the subgraph of $G_i$ consisting of all
vertices and positive edges of $G_i$.
When $G_i$ is of type $\mathcal{S}$, it is obvious that each component of $G_i^+$ has a disk support, that is,
there is a disk in $\widehat{F}_i$ which contains the component in its interior.
Even if $G_i$ is of type $\mathcal{P}$, the same is true, because any orientation-preserving loop in
a projective plane is contractible.

%The following proposition shows
%how the existence of $x$-faces on $G_2$ will be applied.

\begin{proposition} \label{prop:S}
If $G_1$ is of type $\mathcal{S}$ or $\mathcal{P}$
and $G_2$ is any of the six types,
then each non-$sl$-vertex of $G_1$ has at least
$(\Delta(\gamma_1,\gamma_2)-1)n_2+\chi(\widehat{F}_2)$ positive edge endpoints.
\end{proposition}

\begin{proof}
Assume that there is a non-$sl$-vertex $x$ of $G_1$
which has more than $n_2-\chi(\widehat{F}_2)$ negative edges.
Then $G_2$ contains more than $n_2-\chi(\widehat{F}_2)$ positive $x$-edges
by the parity rule.
Thus the subgraph $\Gamma_x$ of $G_2$ consisting of all vertices
and positive $x$-edges of $G_2$ has $n_2$ vertices and
more than $n_2-\chi(\widehat{F}_2)$ edges.
Then an Euler characteristic calculation
shows that $\Gamma_x$ contains a disk face, which is an $x$-face.
This contradicts Theorem \ref{thm:Pxface} or \ref{thm:Sxface}.
\end{proof}

Suppose that $G_1$ is of type $\mathcal{S}$ or $\mathcal{P}$, and that $\Delta(\gamma_1,\gamma_2)\ge 2$.
By Lemmas \ref{lem:n} and \ref{lem:Schar}, $G_1$ has a non-$sl$-vertex.
Note that each of these vertices has therefore
valency at least two in $\overline{G}_1$ by Lemmas \ref{lem:parallel} and \ref{lem:parallel-negative}.
Now take an innermost component $\Lambda_0$ of $G_1^+$
with a disk support $D_0$, which means that $D_0 \cap G_1^+ = \Lambda_0$.
The subgraph $\Lambda_0$ has therefore at most one $sl$-vertex.

Suppose that $\Lambda_0$ is a single vertex.
Then only negative edges are incident there.
By Proposition \ref{prop:S}, it is an $sl$-vertex.
If $G_1$ is of type $\mathcal{S}$,
the edges of a Scharlemann cycle of $G_2$ are incident there.
Then one of the disks bounded by these edges contains non-$sl$-vertices, since $n_1\ge 3$ by Lemma \ref{lem:n}.
Hence we can choose another innermost component of $G_1^+$ with more than one vertex and
at most one $sl$-vertex.
If $G_1$ is of type $\mathcal{P}$, a negative loop is incident there.
Thus we can also choose another innermost component of $G_1^+$
which contains only non-$sl$-vertices.

We may therefore assume that $\Lambda_0$ has more than one vertex.
Then $\Lambda_0$ has no cut vertex or
at least two blocks with at most one cut vertex.
Thus we can choose an innermost component $\Lambda$
with a disk support $D$ after splitting $\Lambda_0$ at all cut vertices,
such that $\Lambda$ has more than one vertex with at most one cut vertex.
In such a case that $\Lambda$ contains a cut vertex and a distinct $sl$-vertex,
we can choose another innermost one containing
no $sl$-vertex.

Such a subgraph $\Lambda$ of $G_1^+$ is called an {\em extremal block\/}
with a disk support $D$.
A vertex of $\Lambda$ is called a {\em ghost vertex\/} if it is either a cut vertex or an $sl$-vertex.
We emphasize that $\Lambda$ has more than one vertex and at most one ghost vertex $y_0$.
A vertex of $\Lambda$ is called a {\em boundary vertex\/} if there is an arc connecting it to $\partial D$ whose interior is
disjoint from $\Lambda$, and an {\em interior vertex\/} otherwise.
Then the preceding argument with Proposition \ref{prop:S}
proves the following theorem which plays key role in this paper:

\begin{theorem} \label{thm:extremalS}
Suppose that $G_1$ is of type $\mathcal{S}$ or $\mathcal{P}$.

\begin{itemize}
\item[(1)] If $G_2$ is of type $\mathcal{S}$, $\mathcal{P}$,
$\mathcal{A}$ or $\mathcal{B}$
with the assumption that $\Delta(\gamma_1,\gamma_2) \geq 2$,
then $G_1$ contains an extremal block $\Lambda$ with a disk support $D$
so that each boundary vertex, except $y_0$, has at least
$n_2+\chi(\widehat{F}_2)$ consecutive edge endpoints of $\Lambda$.
\item[(2)] If $G_2$ is of type $\mathcal{T}$ or $\mathcal{K}$
with the assumption that $\Delta(\gamma_1,\gamma_2) \geq 3$ and $n_2\ge 3$ when $G_2$ is of type $\mathcal{T}$,
and $n_2\ge 2$ when $G_2$ is of type $\mathcal{K}$,
then $G_1$ contains an extremal block $\Lambda$ with a disk support $D$
so that each boundary vertex, except $y_0$,
has at least $2 n_2$ consecutive edge endpoints of $\Lambda$.
\end{itemize}
\end{theorem}

%%%%%%%%%%%%%%%%%%%%%%%%%%%%%%%%%%%%%%%%%%%%%%%%%%%%%%%%%%%%%%%%%
\section{$(\mathcal{S},\mathcal{S})$ case} \label{sec:SS}

\begin{proof}[Proof of Theorem \ref{thm:SS}]
Assume for contradiction that $\Delta(\gamma_1,\gamma_2) \geq 2$.
Theorem \ref{thm:extremalS}(1) says that
$G_1^+$ contains such an extremal block $\Lambda$ with a disk support
that each boundary vertex, except $y_0$,
has more than $n_2$ consecutive edge endpoints,
so different all $n_2$ labels.
We can choose a non-$sl$-label $x$ (at $y_0$ if it exists),
unless $G_2$ is of type $\mathcal{S}$ and
there is $y_0$ with only two edge endpoints which have the $sl$-labels.
Since $\Lambda$ has at least the same number of $x$-edges
as that of vertices, it contains an $x$-face, contradicting Theorem
\ref{thm:Pxface} or \ref{thm:Sxface}.

For the exceptional case, consider $\Lambda-y_0$.
Since all labels still appear on each vertex of $\Lambda-y_0$,
the same argument above leads to a contradiction.
\end{proof}

%%%%%%%%%%%%%%%%%%%%%%%%%%%%%%%%%%%%%%%%%%%%%%%%%%%%%%%%%%%%%%%%%%%%%%%
\section{$(\mathcal{S},\mathcal{A})$ case} \label{sec:SA}

\begin{figure}[ht]
\epsfbox{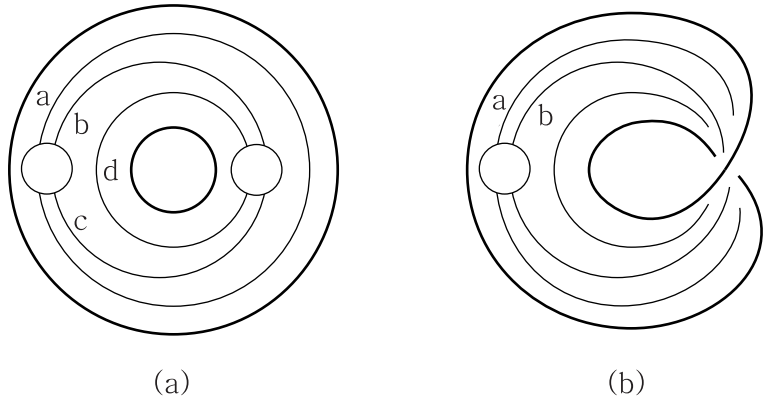}
\caption{Annulus and M\"{o}bius band}\label{fig:annulus}
\end{figure}

\begin{proof}[Proof of Theorem \ref{thm:SA}]
Assume that $\Delta=\Delta(\gamma_1,\gamma_2) \geq 2$.
Recall that $M(\gamma_2)$ is irreducible and
boundary irreducible.

First, assume that $n_2 \geq 3$ when $G_2$ is of type $\mathcal{A}$
(or $n_2 \geq 2$ when of type $\mathcal{B}$).
Theorem \ref{thm:extremalS}(1) says that
$G_1^+$ contains an extremal block $\Lambda$ with a disk support
such that each boundary vertex, except $y_0$, has all different $n_2$ labels.
When $G_2$ is of type $\mathcal{B}$, 
choose a non-$sl$-label $x$ at $y_0$, if it exists.
Otherwise, $x$ is any non-$sl$-label.
Then $\Lambda$ has at least the same number of $x$-edges
as that of vertices, and therefore $\Lambda$ has an $x$-face,
contradicting Theorem \ref{thm:Bxface}.
So $G_2$ is of type $\mathcal{A}$.

\begin{claim}\label{claim:12Schar}
$G_1$ contains a Scharlemann cycle.
\end{claim}

\begin{proof}[Proof of Claim \ref{claim:12Schar}]
Assume not.  Choose any label $x$ (at $y_0$ if it exists).
Then there is an $x$-edge at any vertex of $\Lambda$.
Hence $\Lambda$ contains a great $x$-cycle, and so a Scharlemann cycle \cite[Lemma 2.6.2]{CGLS}.
\end{proof}

Thus we may assume that $G_1$ contains only $12$-Scharlemann cycles.
Equivalently, $sl$-labels of $G_1$ are $1$ and $2$.
By Lemma \ref{lem:Scharsep}, $\widehat{F}_2$ is separating and $n_2$ is even.

\begin{claim}\label{claim:extrem}
\begin{itemize}
\item[(1)] $\Lambda$ has no interior vertices;
\item[(2)] $\Lambda$ has a ghost vertex $y_0$ such that only two edges are incident to $y_0$ in $\Lambda$
with $sl$-labels $1$ and $2$; and
\item[(3)] $n_2=4$, and the two edges incident to $y_0$ are indeed a $14$-edge and a $23$-edge.
\end{itemize}
\end{claim}

\begin{proof}[Proof of Claim \ref{claim:extrem}]
If $\Lambda$ has an interior vertex or has no ghost vertex, then let $x$ be any non-$sl$-label.
Then $\Lambda$ has at least the same number of $x$-edges as that of vertices, so contains an $x$-face, contradicting
Theorem \ref{thm:Axface}.
Even if there is a ghost vertex $y_0$, we can choose a non-$sl$-label $x$ at $y_0$, unless (2) holds.
Finally, assume $n_2\ge 6$.
Then we can choose a non-$sl$-label $x$ such that $y_0$ is not incident to an $x$-edge.
Thus $\Lambda-y_0$ contains an $x$-face as above again.
Even if $n_2=4$, there is still such a label $x$, unless the second conclusion of (3) holds.
\end{proof}

In fact, $\Lambda-y_0$ contains a $12$-Scharlemann cycle $\sigma$, because 
there is a $1$-edge at any vertex of $\Lambda$ by
(3), and hence $\Lambda-y_0$ contains a great $1$-cycle.
By Lemma \ref{lem:Scharsep},
the edges of $\sigma$ cuts $\widehat{F}_2$ into two annuli $A_1,A_2$ and some disks.
%Let $A_i'=\mathrm{cl}\,(\widehat{F}_2-A_i)$ for $i=1,2$.

\begin{claim}\label{claim:34edge}
$\Lambda$ contains a $34$-edge.
\end{claim}

\begin{proof}[Proof of Claim \ref{claim:34edge}]
Assume not.
Then $\Lambda$ contains only $12$, $14$ and $23$-edges by the parity rule.
Let $y_0,y_1,\dots,y_k$ be the vertices of $\Lambda$, numbered consecutively along $\partial \Lambda$, where 
the $14$-edge $e_0$ at $y_0$ is incident to $y_1$.
(Since $y_i$ $(i\ne 0)$ is incident to at least four edges in $\Lambda$, we see $k>1$.)
Let $e_1$ be the $3$-edge at $y_1$, just before $e_0$.
(This means that the two endpoints of $e_0$ and $e_1$ are successive around $y_1$.) 
If $e_1$ goes to $y_k$, then $e_1$ is a $34$-edge.
So assume $e_1$ goes to $y_i$ for $1<i<k$.
Thus $e_1$ has the label $2$ at $y_i$.
Split $\Lambda$ along $e_1$, and let $\Lambda'$ be the part containing $y_1,y_2,\dots,y_i$.
Then $\Lambda'$ has a $3$-face, contradicting Theorem \ref{thm:Axface}.
\end{proof}

Hence the vertices $3$ and $4$ lie in the same component, and
we can see that they lie in $A_1$ or $A_2$ indeed.
Otherwise, the vertex $3$ or $4$ is incident to more than $n_1$ negative edges.
Thus, $G_1$ contains a $3$-face or $4$-face, a contradiction.
We may assume that both are in $A_2$.
Let $A_1'=\mathrm{cl}(\widehat{F}_2-A_2)$, and 
$Y=\mathrm{nhd}(A_1'\cup V_{12}\cup E)$, where $E$ is the face bounded by $\sigma$.
Then the frontier of $Y$ gives an essential annulus (see \cite[Claim on page 430]{W3}), which
meets $V_{\gamma_2}$ in two disks.
This contradicts the minimality of $n_2$.

%%%%%%%%%%%%%%%%%%%%%%%%
If $n_2 = 1$ when $G_2$ is of type $\mathcal{A}$,
$G_2$ has only positive edges, and these edges are all parallel.
This contradicts Lemma \ref{lem:parallel}(1-2).

%%%%%%%%%%%%%%%%%%%%%%%
Finally assume $n_2 = 2$ when $G_2$ is of type $\mathcal{A}$
(or $n_2 = 1$ when of type $\mathcal{B}$).
Assume for contradiction that $\Delta \geq 3$.
Consider the reduced graph $\overline{G}_2$.
An Euler characteristic calculation shows that
each vertex of $\overline{G}_2$ has valency at most 4.
By Lemmas \ref{lem:parallel}(1-2) and  and \ref{lem:parallel-negative}(1) the valency must be 4,
so the graph looks like that in Figure \ref{fig:annulus}
((a) for type $\mathcal{A}$ and (b) for type $\mathcal{B}$).

\begin{claim}\label{claim:annulus2}
In Figure \ref{fig:annulus}(a), the two vertices have opposite signs.
\end{claim}

\begin{proof}[Proof of Claim \ref{claim:annulus2}]
Assume not.  Then $G_2$ has only positive edges.
For a non-$sl$-label $x$ of $G_2$,
let $\Gamma$ be the subgraph of $G_2$ consisting of all vertices and all $x$-edges of $G_2$.
Then an Euler characteristic calculation shows that $\Gamma$ has a disk face, which
is an $x$-face.
This contradicts Theorem \ref{thm:Pxface} or \ref{thm:Sxface}.
\end{proof}

Thus $a$ and $d$ are the families of positive edges
and $b$ and $c$ are the families of negative edges.
When $G_1$ is of type $\mathcal{P}$,
this is impossible by Lemmas \ref{lem:parallel} and \ref{lem:parallel-negative}.
When $G_1$ is of type $\mathcal{S}$,
$G_2$ contains an $S$-cycle by Lemma \ref{lem:parallel},
so $M(\gamma_1)$ contains a projective plane by Lemma \ref{lem:Schar}(3),
a contradiction to the previous result.
\end{proof}

We now give an example realizing the case $\Delta = 2$ in Theorem \ref{thm:SA}.
Theorem 2.6 of \cite{EW} shows that there is a hyperbolic manifold $M$
such that $M(\gamma_1) = (S^1\times D^2)\sharp L(2,1)$ and
$M(\gamma_2)$ is the union of $C(2,1)$ and $Q(2p,-2p)$,
for any integer $p \geq 2$, along a torus, with $\Delta(\gamma_1,\gamma_2)=2$.
We denote by $C(r,s)$ the cable space of type $(r,s)$,
and by $Q(r,s)$ the Seifert fibered manifold
with orbifold a disk with two cone points of index $r$ and $s$.
Note that $C(2,1)$ contains a M\"{o}bius band and essential annulus
with the boundaries on the outside torus, hitting the attached
solid torus once and twice, respectively.

%%%%%%%%%%%%%%%%%%%%%%%%%%%%%%%%%%%%%%%%%%%%%%%%%%%%%%%%%%%%%%%
\section{$(\mathcal{S},\mathcal{T})$ case} \label{sec:ST}

\begin{figure}[ht]
\epsfbox{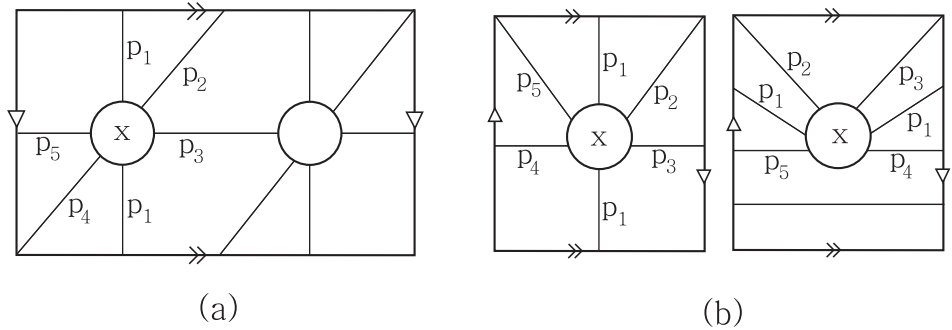}
\caption{Torus and Klein bottle}\label{fig:torus}
\end{figure}

\begin{proof}[Proof of Theorem \ref{thm:ST}]
Assume that $\Delta=\Delta(\gamma_1,\gamma_2)\geq 3$.
Recall that $M(\gamma_2)$ is irreducible
and boundary irreducible.

First, assume that $n_2\geq 3$.
%(Moreover, assume that $n_2\ne 4$ when $G_2$ is of type $\mathcal{T}$.)
Theorem \ref{thm:extremalS}(2) says that
$G_1^+$ contains an extremal block $\Lambda$
with a disk support $D$ so that each boundary vertex, except $y_0$,
has at least $2 n_2$ consecutive edge endpoints of $\Lambda$.
If $G_2$ is of type $\mathcal{K}$, then we can choose a non-$sl$-label $x$ of $G_1$ by Lemma \ref{lem:Schar}.
Assume that $G_2$ is of type $\mathcal{T}$.

\begin{claim}\label{claim:Tcase}
Either $G_1$ has a label $x$, which is not a label of $S$-cycles in $G_1$, or
$M(\gamma_2)$ contains a Klein bottle which meets $V_{\gamma_2}$ at least two times.
\end{claim}

\begin{proof}[Proof of Claim \ref{claim:Tcase}]
If $n_2=3$, then $G_1$ does not have a Scharlemann cycle by Lemma \ref{lem:Scharsep}.
Hence any label is a desired one.
Assume $n_2\ge 4$.
If $G_1$ has four labels of $S$-cycles, then
there are two $S$-cycles whose label pairs are disjoint.
Then $M(\gamma_2)$ contains a Klein bottle $R$ by the argument of the proof of \cite[Lemma 3.10]{GL2}.
If $R$ meets $V_{\gamma_2}$ in a single meridian disk, then take a double covering torus $\widetilde{R}$ of $R$
in $M(\gamma_2)$.
By the assumption $n_2\ge 4$, $\widetilde{R}$ is compressible.
Then $M(\gamma_2)$ is a Seifert fibered manifold, called a prism manifold,
which has finite fundamental group.
This contradicts the fact that $M(\gamma_2)$ is toroidal.
Hence $R$ meets $V_{\gamma_2}$ at least two times.
\end{proof}

In the second conclusion of Claim \ref{claim:Tcase},
if the Klein bottle meets $V_{\gamma_2}$ just two times, then we jump to the case of type $\mathcal{K}$
where $n_2=2$ below.
Thus we have a label $x$ of $G_1$ which is not an $sl$-label when $G_2$ is of type $\mathcal{K}$,
or which is not a label of $S$-cycles otherwise.

%%%%%%%%%%%%%%%%%%%%%%%%%
Let $\Lambda^{x}$ be the subgraph of $\Lambda$
consisting of all vertices and $x$-edges.
Then each boundary vertex of $\Lambda^{x}$, except $y_0$,
has at least two edges attached with label $x$,
which cannot be parallel by Lemma \ref{lem:parallel}(3-4).
Note that $\Lambda^{x}$ may not be connected.
Then, apply the argument in Section \ref{sec:extremal}
to the present situation;
choose an extremal block $\Lambda'$ of $\Lambda^x$ with a disk support $D'$
in $D$, which we can define in a similar way.

Let $v,e$ and $f$ be the numbers of vertices, edges, and disk faces
of $\Lambda'$, respectively.
Also let $v_i, v_{\partial}$ and $v_g$ be the numbers
of interior vertices, boundary vertices and ghost vertices.
Hence $v=v_i+v_{\partial}$ and $v_g=0$ or $1$.

Suppose that $\Lambda'$ has a bigon.
By Lemma \ref{lem:parallel}, it contains either a generalized $S$-cycle or
an extended $S$-cycle.
(Recall that $x$ itself is not a label of $S$-cycle.)
But this is impossible by Lemma \ref{lem:Schar}.
Thus each face of $\Lambda'$ is a disk with at least $3$ sides.
Hence we have $3f + v_{\partial} \leq 2e$.
Since $\Lambda'$ has only disk faces, combined with $v-e+f=\chi(D')=1$,
we get $e \leq 3v_i + 2v_{\partial} - 3$.
On the other hand, we have $2(v_{\partial} - v_g) + \Delta v_i \leq e$,
because each boundary vertex of $\Lambda'$, except $y_0$,
has at least two edges attached with label $x$, and $x$ is not a label of level edges.
These two inequalities give us that $3 \leq 2v_g$, a contradiction.

%%%%%%%%%%%%%%%
Assume that $n_2 = 1$ when $G_2$ is of type $\mathcal{T}$.
Then $G_2$ has only positive edges.
Thus it contains an $x$-face for a non-$sl$-label $x$ of $G_2$,
contradicting Theorems \ref{thm:Pxface} and \ref{thm:Sxface}.

%%%%%%%%%%%%%%%
Assume that $n_2 = 2$ when $G_2$ is of type $\mathcal{K}$.
This case is done by Theorem \ref{thm:extremalS}(2) and the argument of \cite[Section6]{JLOT}, which can be
carried over without change.

%%%%%%%%%%%%%%%%%%%%%%%%%%%%%%%%%%%%%%%%%%%%%%%%%%%%%%%%%%%%%%%
Finally assume that $n_2 = 2$ when $G_2$ is of type $\mathcal{T}$
(or $n_2 = 1$ when of type $\mathcal{K}$).
%The case $(\mathcal{S},\mathcal{T})$ was proved in \cite[Lemma 4.3]{BZ1}, and so
%we consider only the cases $(\mathcal{P},\mathcal{T})$, $(\mathcal{S},\mathcal{K})$ and $(\mathcal{P},\mathcal{K})$.
Assume for contradiction that $\Delta \geq 4$.
An Euler characteristic calculation on $\overline{G}_2$
shows that each vertex has valency at most 6.
Then the graph looks like a subgraph of the graph
shown in Figure \ref{fig:torus}
((a) for type $\mathcal{T}$ and (b) for type $\mathcal{K}$) \cite[Lemma 5.2]{G1}.
Here, $p_i \geq 0$ denotes the number of edges in each parallelism class.

\begin{claim}\label{claim:posi}
When $G_2$ is of type $\mathcal{T}$, all non-loop edges of $G_2$ are negative.
\end{claim}

\begin{proof}[Proof of Claim \ref{claim:posi}]
This follows from the same argument as the proof of Claim \ref{claim:annulus2}.
\end{proof}

Thus $p_1 \leq \frac{n_1}{2}+1$, and $p_i \leq n_1-1$ for $i=2,3,4,5$
by Lemmas \ref{lem:parallel} and \ref{lem:parallel-negative}.
Since $\Delta n_1 \leq (n_1+2)+4(n_1-1)=5n_1-2$,
$\Delta=4$ and all $p_i$'s are non-zero.
Without loss of generality, we can assume that $p_1+p_2+p_3 \geq 2n_1$.
Since $2n_1 < p_1+p_2+p_3+1 \leq (\frac{n_1}{2}+1)+2(n_1-1)+1 < 3 n_1$,
we can write that $p_1+p_2+p_3+1=2n_1+r$, where $0<r<n_1$.
Then $p_1=2n_1+r-(p_2+p_3+1) \geq 2n_1+r-2(n_1-1)-1=r+1$.
Hence $1 \leq r \leq p_1-1$.
Thus the loop family corresponding to $p_1$
have two (non-level) edges with the same label $r$, and so the family contains
an $S$-cycle or a generalized $S$-cycle.
This is impossible when $G_1$ is of type $\mathcal{P}$.
When $G_1$ is of type $\mathcal{S}$, $G_2$ contains an $S$-cycle,
so $M(\gamma_1)$ must be of type $\mathcal{P}$,
a contradiction again.
\end{proof}
%%%%%%%%%%%%%%%%%%%%%%%%%%%%%%%%%%%%%%%%%%%%%%%%%%%%%%%%%%%%%%%%%%%%%%%%%%%%%%%%%%%%%%%%%%%%%%%%

Finally we give an example realizing the case $\Delta=3$ in Theorem \ref{thm:ST}.
Let $M$ be the manifold obtained from the Whitehead link exterior
by Dehn filling one component with slope $6$.
Then $M$ is hyperbolic \cite[p.286]{BZ1} and $M(1)=L(2,1)\sharp L(3,1)$.
Since $M$ contains a once-punctured Klein bottle whose boundary has slope $4$,
$M(4)$ contains a Klein bottle hitting the attached solid torus $V$ once.
Also, the boundary torus of a neighborhood of the Klein bottle in $M(4)$ is known to be incompressible.
Clearly, this torus meets $V$ twice.
(If $M(4)$ contains a torus meeting $V$ once, then
such a torus is non-separating in $M(4)$.
This is impossible, because $H_1(M(4))$ is finite.)

%%%%%%%%%%%%%%%%%%%%%%%%%%%%%%%%%%%%%%%%%%%%%%%%%%%%%%%%%%%%%%%%%%%%%%%%%

\end{document}